\documentclass{article}%
\usepackage{amsmath}
\usepackage{amsfonts}
\usepackage{amssymb}
\usepackage{graphicx}
\usepackage[doublespacing]{setspace}%
\providecommand{\U}[1]{\protect\rule{.1in}{.1in}}
\newtheorem{theorem}{Theorem}

\newtheorem{lemma}[theorem]{Lemma}

\begin{document}

\title{Consistency of M estimates for separable nonlinear regression models}
\author{Mar\'{\i}a Victoria Fasano$^{1}$ (virfeather@yahoo.com.ar)
\and Ricardo Maronna$^{1}$ (rmaronna@retina.ar)\\$^{1}$Departamento de Matem\'{a}tica, Facultad de Ciencias Exactas,\\Universidad de La Plata, C.C. 172, La Plata 1900, Argentina. }
\date{~\ }
\maketitle

\begin{abstract}
Consider a nonlinear regression model : $y_{i}=g\left(  \mathbf{x}%
_{i},\mathbf{\theta}\right)  +e_{i},$ $i=1,...,n,$ where the $\mathbf{x}_{i}$
are random predictors $\mathbf{x}_{i}$ and $\mathbf{\theta}$ is the unknown
parameter vector ranging in a set $\mathbf{\Theta\subset}R^{p}\mathbf{.}$ All
known results on the consistency of the least squares estimator and in general
of M estimators assume that either $\mathbf{\Theta}$ is compact or $g$ is
bounded, which excludes frequently employed models such as the
Michaelis-Menten, logistic growth and exponential decay models. In this
article we deal with the so-called \emph{separable} models, where
$p=p_{1}+p_{2},$ $\mathbf{\theta=}\left(  \mathbf{\alpha,\beta}\right)  $ with
$\mathbf{\alpha\in}A\subset R^{p_{1}},$ $\mathbf{\beta\in}B\subset R^{p_{2},}%
$and $g$ has the form $g\left(  \mathbf{x,\theta}\right)  =\mathbf{\beta}%
^{T}\mathbf{h}\left(  \mathbf{x,\alpha}\right)  $ where $\mathbf{h}$ is a
function with values in $R^{p_{2}}.$ We prove the strong consistency of M
estimators under very general assumptions, assuming that $\mathbf{h}$ is a
bounded function of $\mathbf{\alpha,}$ which includes the three models
mentioned above.

Key words and phrases: Nonlinear regression, separable models, consistency,
robust estimation.

\end{abstract}

\section{Introduction}

Consider i.i.d. observations $\left(  \mathbf{x}_{i},y_{i}\right)  ,$
$\ i=1,...,n,$ given by the nonlinear model with random predictors:%
\begin{equation}
y_{i}=g\left(  \mathbf{x}_{i},\mathbf{\theta}_{0}\right)  +e_{i}%
,\ \label{nonlinmod}%
\end{equation}
where $\mathbf{x}_{i}\in R^{q}$ and $e_{i}$ are independent, and the unknown
parameter vector $\mathbf{\theta}_{0}$ ranges in a set $\Theta\subset R^{p}.$
An important case, usually called \emph{separable}, are models where
$p=p_{1}+p_{2}$ and $\mathbf{\theta}_{0}\mathbf{=}\left(  \mathbf{\alpha}%
_{0}\mathbf{,\beta}_{0}\right)  $ with $\mathbf{\alpha}_{0}\mathbf{\in
}A\subset R^{p_{1}}$ and $\mathbf{\beta}_{0}\mathbf{\in}B\subset R^{p_{2}},$
and $g$ of the form%
\begin{equation}
g\left(  \mathbf{x,\theta}\right)  =g\left(  \mathbf{x,\alpha,\beta}\right)
=\sum_{j=1}^{p_{2}}\beta_{j}h_{j}\left(  \mathbf{x,\alpha}\right)  ,
\label{lincompo}%
\end{equation}
where $h_{j}$ $(j=1,...,p_{2})$ are functions of $X\times R^{p_{2}}\rightarrow
R.$ Usually $B$ is the whole of $R^{p_{2}}$ or an unbounded subset of it.
Examples are the Michaelis-Menten model, with%
\begin{equation}
p_{1}=p_{2}=q=1,\ x\geq0,~\alpha,\beta>0,~h\left(  x,\alpha\right)  =\frac
{x}{x+\alpha}, \label{Michael}%
\end{equation}
the logistic growth model, with%
\begin{equation}
q=1,\ p2=1,\ \ p_{2}=1,\ x\geq0,\ \alpha_{j}>0,\ \beta>0,\ h\left(
\mathbf{x,\alpha}\right)  =\frac{e^{\alpha_{2}x}}{1+\alpha_{1}\left(
e^{\alpha_{2}x}-1\right)  }, \label{logis}%
\end{equation}
the exponential decay model, with
\begin{equation}
~q=1,~p_{2}=p_{1}+1,\ ~x\geq0,~\alpha_{j}<0,~\beta_{j}\geq0,~g\left(
\mathbf{x,\alpha,\beta}\right)  =\beta_{0}+\sum_{j=1}^{p_{1}}\beta
_{j}e^{\alpha_{j}x}, \label{expdecay}%
\end{equation}
and the exponential growth model, like (\ref{expdecay}) but with $\alpha
_{j}>0.$

The classical least squares estimate (LSE) is given by%
\[
\widehat{\mathbf{\theta}}=\arg\min_{\mathbf{\theta\in}\Theta}\sum_{i=1}%
^{n}\left(  y_{i}-g\left(  \mathbf{x}_{i},\mathbf{\theta}\right)  \right)
^{2}.
\]

The consistency of the LSE assuming $\mathrm{E}\left(  e_{i}\right)  =0$ and
$\mathrm{Var}\left(  e_{i}\right)  =\sigma^{2}<\infty$ has been proved by
several authors under the assumption of a compact $\Theta$; in particular
Amemiya (1983), Jennrich (1969) and Johansen (1984). Wu (1981) assumes that
$\Theta$ is a finite set.

Richardson and Bhattacharyya (1986) do not require the compactness of
$\Theta,$ but they assume $g\left(  \mathbf{x,\theta}\right)  $ to be a
bounded function of $\mathbf{\theta,}$ which excludes most separable models.

Shao (1992) showed the consistency of the LSE without requiring the compacity
of $\Theta$ nor the boundedness of $g,$ but requires assumptions on $g$ that
exclude the simplest separable models. For example, in the case $g\left(
x,\mathbf{\theta}\right)  =\beta e^{\alpha x},$ for any $x_{0}>0$ one can make
$g\left(  x_{0},\mathbf{\theta}\right)  =$constant with $\alpha\rightarrow
-\infty$ and $\beta\rightarrow0.$ This fact violates both \textquotedblleft
Condition 1\textquotedblright\ and \textquotedblleft Condition
2\textquotedblright\ in page 427 of his paper.

The well-known fact that the LSE is sensitive to outliers has led to the
development of \emph{robust estimates} that are simultaneously highly
efficient for normal errors and resistant to perturbations of the model. One
of the most important families of robust estimates are the \emph{M-estimates}
proposed by Huber (1973) for the linear model. For nonlinear models they are
defined by%
\begin{equation}
\mathbf{\hat{\theta}}_{n}\mathbf{=}\arg\min_{\mathbf{\theta}\in\Theta}%
\sum_{i=1}^{n}\rho\left(  \frac{y_{i}-g\left(  \mathbf{x}_{i},\mathbf{\theta
}\right)  }{\widehat{\sigma}}\right)  , \label{defMestimaSig}%
\end{equation}
where $\rho$ is a loss function whose properties will be described in the next
section and $\widehat{\sigma}$ is an estimate of the error's scale. However,
at this stage of our research we deal with the simpler case of known $\sigma.$
Then it may be assumed without loss of generality that $\sigma=1$ and
therefore we shall deal with estimates of the form%
\begin{equation}
\mathbf{\hat{\theta}}_{n}\mathbf{=}\arg\min_{\mathbf{\theta}\in\Theta}%
\sum_{i=1}^{n}\rho\left(  y_{i}-g\left(  \mathbf{x}_{i},\mathbf{\theta
}\right)  \right)  . \label{defMestima}%
\end{equation}

All published results on the consistency of robust estimates for nonlinear
models require the compacity of $\Theta.$ Oberhofer (1982) deals with the
$L_{1}$ estimator. Vainer and Kukush (1998) and Liese and Vajda (2003, 2004)
deal with M estimates. The latter deal with $O\left(  n^{-1/2}\right)  $
consistency and asymptotic normality of M estimates in more general models.
Stromberg (1995) proved the consistency of the Least Median of Squares
estimate (Rousseeuw, 1984), and \v{C}\'{\i}\v{z}ek (2005) dealt with the
consistency and asymptotic normality of the Least Trimmed Squares estimate.
Fasano et al. (2012) study the functionals related to M estimators in linear
and nonlinear regression; in the latter case, they also assume a compact
$\mathbf{\Theta.}$

In this article we will prove the consistency of M estimates for separable
models without assuming the compactness of $\Theta$, but assuming the
boundedness of the $h_{j}$s$\mathbf{;}$ this case includes the exponential
decay, logistic growth and Michaelis-Menten models. It can thus be considered
as a generalization of (Richardson and Bhattacharyya, 1986).

\section{The assumptions}

It will be henceforth assumed that $\rho$ is a \textquotedblleft$\rho
$--function\textquotedblright\ in the sense of (Maronna et al, 2006). i.e.,
$\rho\left(  u\right)  $ is a continuous nondecreasing function of $|u|,$ such
that $\rho\left(  0\right)  =0$ and that if $\rho(u)<\sup_{u}\rho(u)$ and
$0\leq u<v$ then $\rho(u)<\rho(v).$ We shall consider two cases: unbounded
$\rho$ and bounded $\rho.$ The first includes convex function, in particular
the LSE with $\rho\left(  x\right)  =x^{2}$ and the well-known Huber function%
\begin{equation}
\rho_{k}(x)=\left\{
\begin{array}
[c]{ccc}%
x^{2} & \mathrm{if} & \left\vert x\right\vert \leq k\\
2k\left\vert x\right\vert -k^{2} & \mathrm{if} & \left\vert x\right\vert >k
\end{array}
\right.  \label{dfhubrho}%
\end{equation}
and the second includes the bisquare function $\rho\left(  x\right)
=\min\left\{  1-\left(  1-\left(  x/k\right)  ^{2}\right)  ^{3},1\right\}  $,
where $k$ is in both cases a constant that controls the estimator's efficiency.

Let $\mathbf{h}\left(  \mathbf{x,\alpha}\right)  =\left(  h_{1}\left(
\mathbf{x,\alpha}\right)  ,...,h_{p_{2}}\left(  \mathbf{x,\alpha}\right)
\right)  ^{\prime}$ where in general $\mathbf{a}^{\prime}$ denotes the
transpose of $\mathbf{a.}$The necessary assumptions are:

\begin{description}
\item[A] $B$ is a closed set such that $t\mathbf{\beta\in}B$ for all
$\mathbf{\beta\in}B$ and $t>0.$

\item[B] $\sup_{\mathbf{\alpha\in}A}\mathrm{E}|\rho\left(  y-\mathbf{\beta
}^{\prime}\mathbf{h}\left(  \mathbf{x,\alpha}\right)  \right)  |<\infty$ for
all $\mathbf{\ \beta\in}B.$

\item[C] The function $\mathrm{E}\rho\left(  e-t\right)  $ --where $e$ denotes
any copy of $e_{i}$-- has a unique minimum at $t=0.$ Put $\lambda
_{0}=\mathrm{E}\rho\left(  e\right)  .$

\item[D] $\mathbf{h}$ is continuous in $\mathbf{\alpha}$ a.s. and%
\begin{equation}
\mathbf{\alpha}\not =\mathbf{\alpha}_{0}\Rightarrow\sup_{\mathbf{\beta}\in
B}\mathrm{P}\left\{  \mathbf{\beta}^{\prime}\mathbf{h}\left(  \mathbf{x,\alpha
}\right)  =\mathbf{\beta}_{0}^{\prime}\mathbf{h}\left(  \mathbf{x,\alpha}%
_{0}\right)  \right\}  <1\text{ } \label{CL1}%
\end{equation}

\item[E] Let $S=\sup_{t}\rho\left(  t\right)  $ (which may be infinite). Then%
\begin{equation}
\delta=:\sup_{\mathbf{\beta\not =0,}\text{ }\mathbf{\alpha}\in A}%
\mathbf{P}\left(  \mathbf{\beta}^{\prime}\mathbf{h}\left(  \mathbf{x,\alpha
}\right)  =0\right)  <1-\frac{\lambda_{0}}{S}. \label{CL2}%
\end{equation}

\item[F] Call $\mathcal{U}$ the family of all open neighborhoods of
$\mathbf{\alpha}_{0}.$ Then
\[
\sup_{\mathbf{\beta}}\inf_{U\in\mathcal{U}}\sup_{\mathbf{\alpha\notin}%
U}\mathrm{P}\left\{  \mathbf{\beta}^{\prime}\mathbf{h}\left(  \mathbf{x,\alpha
}\right)  =\mathbf{\beta}_{0}^{\prime}\mathbf{h}\left(  \mathbf{x,\alpha}%
_{0}\right)  \right\}  <1.
\]

\item[G] $\mathbf{h}$ is bounded as a function of $\mathbf{\alpha,}$ i.e.,
$\sup_{\mathbf{\alpha}\in A}\left\Vert \mathbf{h}\left(  \mathbf{x,\alpha
}\right)  \right\Vert <\infty\ \mathrm{a.s.}$
\end{description}

We now comment on the assumptions.

For (A) to hold in examples (\ref{Michael})-(\ref{logis})-(\ref{expdecay}) we
must enlarge the range of $\beta_{j}$s to $\beta_{j}\geq0.$ However, to ensure
the validity of (D) and (F), it will be assumed that the elements of the
\textquotedblleft true\textquotedblright\ vector $\mathbf{\beta}_{0}$ are all positive.

If $\rho$ is bounded, (B) holds without further conditions. Sufficient
conditions for Huber's $\rho$ and for the LSE are finite moments of $e$ and of
$\mathbf{h}\left(  \mathbf{x,\alpha}\right)  ,$ of orders one and two, respectively.

A sufficient condition for (C) is that the distribution of $e$ has an even
density $f\left(  u\right)  $ that is nonincreasing for $u\geq0$ and is
decreasing in a neighborhood of $u=0$ (see Lemma 3.1 of Yohai (1987)).
If\ $\rho$ is strictly convex with a derivative $\psi,$ then a sufficient
condition is $\mathrm{E}\psi\left(  e\right)  =0,$ which for the LSE reduces
to $\mathrm{E}e=0.$

Assumption (D) is required for ensure uniqueness of solutions. For examples
(\ref{Michael})-(\ref{logis}) it is very easy to verify. For (\ref{expdecay})
it follows from the well-known linear independence of exponentials.

If $S=\infty,$ (E) just means that $\delta<1$ (since $\lambda_{0}<\infty$ by
(B)). Otherwise it puts a bound on $\delta.$ In our examples we have
$\delta=0,$ since $\mathbf{\beta}^{\prime}\mathbf{h}>0$ if $\mathbf{\beta}$
has a single nonnull (positive) element.

Assumption (F) is required in the case of non-compact $A,$ to prevent the
estimator $\widehat{\mathbf{\alpha}}$ from \textquotedblleft escaping to the
border\textquotedblright. In our examples the border for the $\alpha_{j}$s is
either zero of infinity, and (F) is easily verified by a detailed but
elementary calculation (taking into account the remark above that all elements
of $\mathbf{\beta}_{0}$ are positive). For example, in (\ref{Michael}) it
suffices to consider neighborhoods of the form $\left(  \alpha_{0}%
/K,K\alpha_{0}\right)  $ with $K$ sufficiently large.

Finally, (G) is easily verified for models (\ref{Michael})-(\ref{logis}%
)-(\ref{expdecay}).

\section{The results}

For separable\emph{ }models\emph{ } the M-estimate is given by
\[
\mathbf{\hat{\theta}}_{n}=\left(  \widehat{\mathbf{\alpha}}_{n},\widehat
{\mathbf{\beta}}_{n}\right)  \mathbf{=}\arg\min_{\mathbf{\alpha\in}A,\text{
}\mathbf{\beta}\in B}\frac{1}{n}\sum_{i=1}^{n}\rho\left(  y_{i}-\mathbf{\beta
}^{\prime}\mathbf{h}\left(  x_{i},\mathbf{\alpha}\right)  \right)  .
\]

We now state our main result.

\begin{theorem}
Assume model (\ref{lincompo}) with conditions A-B-C-D-E-F-G. Then the M
estimate $\left(  \widehat{\mathbf{\alpha}}_{n},\widehat{\mathbf{\beta}}%
_{n}\right)  $ is strongly consistent for $\mathbf{\theta}_{0}.$
\end{theorem}

We shall first need an auxiliary result, based on a proof in (Bianco and
Yohai, 1996).

\begin{lemma}
Assume model (\ref{lincompo}) with conditions A-B-C-D-E and $A$ \emph{compact}%
. Then $\left\Vert \mathbf{\hat{\beta}}_{n}\right\Vert $ is ultimately bounded
with probability one.
\end{lemma}

\textbf{Proof of the Lemma: }Put%
\[
\lambda\left(  \mathbf{\alpha,\beta}\right)  =\mathrm{E}\rho\left(
y-\mathbf{\beta}^{\prime}\mathbf{h}\left(  \mathbf{x,\alpha}\right)  \right)
.
\]

It follows from (C) that $\mathbf{\lambda(\alpha},\mathbf{\beta)}$ attains its
minimum only when $\mathbf{\beta}^{\prime}\mathbf{h}\left(  \mathbf{x,\alpha
}\right)  =\mathbf{\beta}_{0}^{\prime}\mathbf{h}\left(  \mathbf{x,\alpha}%
_{0}\right)  $ a.s. and by (\ref{CL1}) this happens when $\left(
\mathbf{\alpha},\mathbf{\beta}\right)  =\left(  \mathbf{\alpha}_{0}%
,\mathbf{\beta}_{0}\right)  .$ Therefore
\begin{equation}
\left(  \mathbf{\alpha},\mathbf{\beta}\right)  \neq\left(  \mathbf{\alpha}%
_{0},\mathbf{\beta}_{0}\right)  \Rightarrow\lambda\left(  \mathbf{\alpha
},\mathbf{\beta}\right)  >\lambda\left(  \mathbf{\alpha}_{0},\mathbf{\beta
}_{0}\right)  =\lambda_{0}. \label{unique}%
\end{equation}

Let $\Gamma=\left\{  \mathbf{\gamma}\in B:\left\Vert \mathbf{\gamma
}\right\Vert =1\right\}  .$ Then we may write $\mathbf{\beta}=t\mathbf{\gamma
}$ with $t=\left\Vert \mathbf{\beta}\right\Vert \in R_{+}$ and $\mathbf{\gamma
}\in\Gamma.$

We divide the proof into two cases.

\textbf{Case I:\ bounded }$\rho:$ Assume that $S=\sup_{u}\rho\left(  u\right)
<\infty.$ To simplify notation it will be assumed without loss of generality
that $S=1.$ For each $\left(  \mathbf{\alpha},\mathbf{\gamma}\right)  \in
A\times\Gamma$ we have%
\[
\lim_{t\rightarrow\infty}\mathrm{E}\rho\left(  \mathbf{y-}t\mathbf{\gamma
}^{\prime}\mathbf{h}\left(  \mathbf{x,\alpha}\right)  \right)  \geq
1-\delta>\lambda_{0},
\]
where $\delta$ is defined in (\ref{CL2}). Let
\[
\xi=1-\delta-\lambda_{0}>0,\ \ \varepsilon=\frac{\xi}{4}<\frac{1-\delta}{4}.
\]

Since (\ref{CL2}) implies that $\mathrm{P}\left(  \left\vert \mathbf{\gamma
}^{\prime}\mathbf{h}\left(  \mathbf{x,\alpha}\right)  \right\vert >0\right)
\geq1-\delta$ for $\mathbf{\gamma\in}\Gamma,$ then for each $\left(
\mathbf{\alpha},\mathbf{\gamma}\right)  \in A\times\Gamma$ there are positive
$a,b$ such that%
\begin{equation}
\mathrm{P}\left(  \left\vert y\right\vert \leq a,\left\vert \mathbf{\gamma
}^{\prime}\mathbf{h}\left(  \mathbf{x,\alpha}\right)  \right\vert \geq
b\right)  \geq1-\delta-\varepsilon. \label{CL3}%
\end{equation}

Then by (\ref{CL3}) there exists $T>0$ such that $t>T$ implies%
\begin{equation}
\mathrm{E}\inf_{t>T}\rho\left(  \mathbf{y-}t\mathbf{\gamma}^{\prime}%
\mathbf{h}\left(  \mathbf{x,\alpha}\right)  \right)  >1-\delta-2\varepsilon.
\label{CL4}%
\end{equation}

Therefore (\ref{CL4}) implies that for each $\left(  \mathbf{\alpha
},\mathbf{\gamma}\right)  \in A\times\Gamma$ there exist a neighborhood
$U\left(  \mathbf{\alpha},\mathbf{\gamma}\right)  \subset A\times\Gamma$ and
$T\left(  \mathbf{\alpha},\mathbf{\gamma}\right)  \in R_{+}$ such that
\begin{equation}
\mathrm{E}\inf_{\left(  \mathbf{\alpha}_{1}\mathbf{,\gamma}_{1}\right)  \in
U\left(  \mathbf{\alpha,\gamma}\right)  }\inf_{\ t>T\left(  \mathbf{\alpha
,\gamma}\right)  }\rho\left(  \mathbf{y-}t\mathbf{\gamma}_{1}^{\prime
}\mathbf{h}\left(  \mathbf{x,\alpha}_{1}\right)  \right)  >1-\delta
-2\varepsilon=\lambda_{0}+\frac{\xi}{2}. \label{CL5}%
\end{equation}

The neighborhoods $\{U\left(  \mathbf{\alpha,\gamma}\right)  :\mathbf{\alpha
\in}A,\mathbf{\gamma}\in\Gamma\}$ are a covering of the compact set
$A\times\Gamma,$ and therefore there exists a finite subcovering thereof:
$\left\{  U_{j}=U\left(  \alpha_{j},\mathbf{\gamma}_{j}\right)  \right\}
_{j=1}^{N}$. Let $T_{0}=\max_{j}T\left(  \alpha_{j},\mathbf{\gamma}%
_{j}\right)  .$

We shall show that $\lim\sup_{n\rightarrow\infty}\left\Vert \mathbf{\hat
{\beta}}_{n}\right\Vert \leq T_{0}$ a.s. Put for brevity%
\[
\lambda_{n}\left(  \mathbf{\alpha,\beta}\right)  =\frac{1}{n}\sum_{i=1}%
^{n}\rho\left(  y_{i}-\mathbf{\beta}^{\prime}\mathbf{h}\left(  x_{i}%
,\mathbf{\alpha}\right)  \right)  .
\]
Then%
\begin{align*}
\inf_{\left\Vert \mathbf{\beta}\right\Vert >T_{0}}\inf_{\mathbf{\alpha}\in
A}\lambda_{n}\left(  \mathbf{\alpha,\beta}\right)   &  \geq\frac{1}{n}%
\sum_{i=1}^{n}\inf_{\mathbf{\alpha}\in A,\mathbf{\gamma}\in\Gamma}%
\inf_{t>T_{0}}\rho\left(  y_{i}-t\mathbf{\gamma}^{\prime}\mathbf{h}\left(
x_{i},\mathbf{\alpha}\right)  \right) \\
&  =\min_{j=1,...,N}\frac{1}{n}\sum_{i=1}^{n}\inf_{\left(  \mathbf{\alpha
,\gamma}\right)  \in U_{j}}\inf_{t>T_{0}}\rho\left(  y_{i}-t\mathbf{\gamma
}^{\prime}\mathbf{h}\left(  x_{i},\mathbf{\alpha}\right)  \right)  ,
\end{align*}
and therefore (\ref{CL5}) and the Law of Large Numbers imply%
\[
\lim\inf_{n\rightarrow\infty}\inf_{\left\Vert \mathbf{\beta}\right\Vert
>T_{0}}\inf_{\mathbf{\alpha}\in A}\lambda_{n}\left(  \mathbf{\alpha,\beta
}\right)  \geq\lambda_{0}+\frac{\xi}{2}~\mathrm{a.s.},
\]
while%
\[
\lambda_{n}\left(  \widehat{\mathbf{\alpha}}_{n}\mathbf{,}\widehat
{\mathbf{\beta}}_{n}\right)  =\inf_{\beta\in B}\inf_{\mathbf{\alpha}\in
A}\lambda_{n}\left(  \mathbf{\alpha,\beta}\right)  \leq\lambda_{n}\left(
\mathbf{\alpha}_{0}\mathbf{,\beta}_{0}\right)  \rightarrow\lambda
_{0}~\mathrm{a.s.}%
\]
which shows that ultimately $\left\Vert \mathbf{\hat{\beta}}_{n}\right\Vert
\leq T_{0}$ with probability one.

\textbf{Case II:\ unbounded }$\rho:$ Here an analogous but simpler procedure
shows the existence of $T_{0}$ and neighborhoods $U\left(  \mathbf{\alpha
,\gamma}\right)  $ such that the left-hand member of (\ref{CL5}) is larger
than $2\lambda_{0},$ and the rest of the proof is similar.$\blacksquare$

\textbf{Proof of the Theorem: }If $A$ is not compact, we employ the same
approach as in (Richardson and Bhattacharyya, 1986): the \v{C}ech-Stone
compactification yields a compact set $\widetilde{A}\supset A$ such that each
bounded continuous function on $A$ has a unique continuous extension to
$\widetilde{A}.$ We have to ensure that (B), (D) and (E) continue to hold for
$\mathbf{\alpha\in}\widetilde{A}.$ Since each element of $\widetilde{A}$ is
the limit of a sequence of elements of $A,$ (B) and (E) are immediate; and (D)
follows from assumption (F). Therefore we can apply the Lemma to conclude that
$\left(  \widehat{\mathbf{\alpha}}_{n}\mathbf{,}\widehat{\mathbf{\beta}}%
_{n}\right)  $ remains ultimately in a compact a.s. The Theorem then follows
from Theorem 1 of Huber (1967).$\blacksquare$

\section{Acknowledgements:}

\textbf{ }This research was partially supported by grants PID 5505 from
CONICET and PICTs 21407 and 00899 from ANPCYT, Argentina.

\medskip

\textbf{References}

Bianco, A., Yohai, V.J., 1996. Robust estimation in the logistic re\-gression
model, in Robust Statistics, Data Analysis and Computer Intensive Methods,
Proceedings of the workshop in honor of Peter J. Huber, editor H. Rieder,
Lecture Notes in Statistics 109, 17-34 Springer-Verlag, New York.

\v{C}\'{\i}\v{z}ek, P., 2006. Least trimmed squares in nonlinear regression
under dependence. Jr. Statist. Plann. \& Inf., 136, 3967-3988.

Fasano, M.V., Maronna, R.A., Sued, M., Yohai, V.J., 2012. Continuity and
differentiability of regression M functionals. Bernouilli (to appear).

Huber, P. J., 1967. The behavior of maximum likelihood estimates under
nonstandard conditions, in Proceedings of the Fifth Berkeley Symposium in
Mathematical Statistics and Probability, Berkeley: University of California
Press, Vol. 1, 221-233.

Jennrich, R. I., 1969. Asymptotic properties of nonlinear least squares
estimators. Ann. Math. Statist., 40\textbf{,} 633-643.

Liese, F., Vajda, I., 2003. A general asymptotic theory of M-estimators I.
Math. Meth. Statist., 12, 454-477.

Liese, F. Vajda, I., 2004. A general asymptotic theory of M-estimators II.
Math. Meth. Statist., 13, 82-95.

Maronna, R.A., Martin, R.D., Yohai, V.J., 2006. Robust Statistics: Theory and
Methods, John Wiley and Sons, New York.

Oberhofer, W., 1982. The consistency of nonlinear regression minimizing the
$L_{1}$ norm. Ann. Statist., 10, 316-319.

Richardson, G.D., Bhattacharyya, B.B., 1986. Consistent estimators in
nonlinear regression for a noncompact parameter space. Ann. Statist.,
14,\textbf{ }1591-1596.

Rousseeuw, P., 1984. Least median of squares regression. Jr.Amer. Statist.
Assoc., 79, 871-880.

Shao, J., 1992. Consistency of Least-Squares Estimator and Its Jackknife
Variance Estimator in Nonlinear Models. Can. Jr. Statist., 20, 415-428.

Stromberg, A. J., 1995. Consistency of the least median of squares estimator
in nonlinear regression. Commun. Statist.: Th. \& Meth., 24, 1971-1984.

Tabatabai M. A.,Argyros I. K., 1993. Robust estimation and testing for general
nonlinear regression models. Appl. Math. \& Comp., 58, 85-101.

Vainer, B. P., Kukush, A. G., 1998. The consistency of M-estimators
constructed from a concave weight function. Th. Prob. \& Math. Statist., 57, 11-18.

Wu, C. F., 1981. Asymptotic theory of nonlinear least squares estimation. Ann.
Statist., 9, 501-513.

Yohai, V. J., 1987. High breakdown-point and high efficiency estimates for
regression. Ann. Statist., 15, 642-656.
\end{document}